\documentclass{amsart}

\newcommand{\Z}{\mathbb{Z}}

\def\calO{\mathcal{O}}

\newtheorem{thm} {Theorem}

\newtheorem{rem} [thm]{Remark}

\begin{document}

\title[]{Primes in the denominators\\ of Igusa class polynomials}

\author{Kristin E. Lauter}
\address{Microsoft Research,
         One Microsoft Way,
         Redmond, WA 98052, USA.}
\email{klauter@microsoft.com}

\date{November 27, 2002}


\maketitle
\section{Introduction}
The purpose of this note is to suggest an analogue for genus $2$ curves of part of Gross and 
Zagier's work on elliptic curves \cite{GZ}.  Experimentally, for genus $2$ curves with CM by 
a quartic CM field $K$, it appears that primes dividing 
the denominators of the discriminants of the Igusa class polynomials all have the property 1) 
that they are bounded by $d$, the absolute value of the discriminant of $K$, and 
2) that they divide $d-x^2$, for some integer $x$ whose square is less than $d$.  A slightly
stronger condition is given in Section 3.
Such primes are primes of bad reduction for the genus $2$ curve and primes of supersingular 
reduction for the Jacobian of the genus $2$ curve.

\section{Algorithm to generate Igusa class polynomials of quartic CM fields}
\label{algo}

\noindent
The steps of this algorithm are as follows:

\begin{enumerate}
\item Choose a CM field $K$ of degree $4$ which is either not Galois or Galois with
Galois group $\Z/4\Z$. Let $K_0$ be the totally real quadratic subfield. 
Choose a CM-type, $\Phi$ (a choice of two complex embeddings of $K$, such that one is not
the complex conjugate of the other).

\item For each element of the ideal class group of $\calO_K$,
choose a representative ideal $\mathfrak{a}_i$, and find an integral
basis $\{1, \tau_i\}$ for it over the ring of integers of $K_0$.

\item Each ideal class corresponds to an isomorphism class of an
abelian variety via Shimura's theory \cite[p.126]{Shimura}.  
The principal polarization(s) on the abelian variety and
the corresponding $2 \times 2$ period matrix (matrices) $\Omega_i$ are given in
\cite[Section 4.2 and p.62]{Spallek}.
The entries of $\Omega_i$ are given in terms of
$\tau_i$, the CM-type $\Phi$, and $\omega$, the generator of the ring of integers of $K_0$.

\item For each period matrix, evaluate the ten even theta constants up to some amount 
of precision. 

\item The three absolute Igusa invariants, $j_1$, $j_2$, $j_3$, associated to each period matrix 
are defined as combinations of the ten even theta constants by formulas given in the appendix.

\item Take the product over all possible period matrices to form the three 
Igusa class polynomials:
$$h_1(x) = \prod_{\Omega_i}(x-j_1(\Omega_i)),$$
$$h_2(x) = \prod_{\Omega_i}(x-j_2(\Omega_i)),$$
$$h_3(x) = \prod_{\Omega_i}(x-j_3(\Omega_i)).$$

\end{enumerate}

This algorithm was implemented using the software packages Pari and MAPLE in \cite{CL}.
Implementations can be found in the literature in \cite{Spallek}, \cite{vW}, \cite{Weng}.

\begin{rem}
The above algorithm is analogous to computing the Hilbert class polynomial associated to
an imaginary quadratic field $K$.  The $j$-invariant for each ideal class is
evaluated up to some amount of precision, then $H(x)$ is formed
by taking the product over all ideal classes of $(x-j(\tau_i))$, where $\tau_i$
is a particular algebraic integer associated to each ideal class.
Then the Hilbert class polynomial, $H(x)$, has integer coefficients.

The case of a quartic CM field is different because there are three class polynomials, 
their coefficients are rational numbers, and the amount of precision required to compute them is 
not known in advance.  The next section proposes a constraint on the primes appearing in the 
denominators which, together with a bound on the power to which each prime appears, would 
give a bound on the amount of precision required to compute the coefficients of the class 
polynomials.
\end{rem}

\section{Conjectural formula for the primes in the denominators of the Igusa 
class polynomials}

Let $h_1$, $h_2$, and $h_3$ be the three Igusa class polynomials with rational coefficients 
obtained from the algorithm in Section \ref{algo}.
Let $\{q_i\}$ be the primes appearing in the denominators of the discriminants 
of all three $h_i$. Let $d$ be the absolute value of the 
discriminant of the number field $K$ and 
let $d_0$ be the absolute value of the norm of the discriminant of $K/K_0$.

Experimentally I have observed that the primes $q_i$ satisfy the following property:

\vskip .1 truein
\noindent
{\bf Property (1)} Each $q_i$ divides $d-x^2$, for some $x$, an integer such that $x^2 \le d$.
In fact in every case so far, each $q_i$ divides $d_0-x^2$, for some $x$, 
an integer such that $x^2 \le d_0$.

\noindent




\begin{rem}
I have only tested this property on a few handfuls of examples, but Annegret Weng has also 
tested it on at least that many others.  
I am currently working on a proof of this property jointly 
with Annegret Weng, Farshid Hajir, Fernando Rodriguez-Villegas, and Tonghai Yang. 

\end{rem}

\section{Appendix: Igusa invariants}

Let $\{\theta_i\}_{i=1,10}$ be the ten even theta constants associated to a given period matrix.
Define the functions $f$ and $g$ as follows:
$$f(k_1,k_2,k_3,k_4,k_5,k_6,k_7,k_8) = (\prod_{i=1,...,8} \theta_{k_i})^4,$$

$$g(k_1,k_2,k_3,k_4,k_5,k_6) = (\prod_{i=1,...,6} \theta_{k_i})^4.$$
Now $h_4$, $h_{10}$, $h_{12}$, and $h_{16}$ are values of modular
forms of weights $4$, $10$, $12$, and $16$:

$h_4 = \sum_{i=1,...,10}\theta_i^8$,

$h_{10} = \prod_{i=1,...,10}\theta_i^2$,

$h_{12} = g(1,5,2,9,6,10)+g(1,2,9,6,8,3)+g(5,9,6,8,10,7)+
g(5,2,6,8,3,7)+g(1,5,2,10,3,7)+g(1,9,8,10,3,7)+g(1,5,2,8,10,4)+
g(1,5,9,8,3,4)+g(5,9,6,10,3,4)+g(2,6,8,10,3,4)+g(1,2,9,6,7,4)+
g(1,5,6,8,7,4)+g(2,9,8,10,7,4)+g(5,2,9,3,7,4)+g(1,6,10,3,7,4),$

$h_{16} = f(8,1,5,2,9,6,8,10)+f(5,1,5,2,9,6,8,3)+f(10,1,2,9,6,8,10,3)+
f(3,1,5,2,9,6,10,3)+f(1,1,5,9,6,8,10,7)+f(2,5,2,9,6,8,10,7)+
f(1,1,5,2,6,8,3,7)+f(9,5,2,9,6,8,3,7)+f(9,1,5,2,9,10,3,7)+
f(6,1,5,2,6,10,3,7)+f(5,1,5,9,8,10,3,7)+f(2,1,2,9,8,10,3,7)+
f(6,1,9,6,8,10,3,7)+f(8,1,5,2,8,10,3,7)+f(10,5,2,6,8,10,3,7)+
f(3,5,9,6,8,10,3,7)+f(7,1,5,2,9,6,10,7)+f(7,1,2,9,6,8,3,7)+
f(9,1,5,2,9,8,10,4)+f(6,1,5,2,6,8,10,4)+f(2,1,5,2,9,8,3,4)+
f(6,1,5,9,6,8,3,4)+f(1,1,5,9,6,10,3,4)+f(2,5,2,9,6,10,3,4)+
f(1,1,2,6,8,10,3,4)+f(5,5,2,6,8,10,3,4)+f(9,2,9,6,8,10,3,4)+
f(8,5,9,6,8,10,3,4)+f(10,1,5,9,8,10,3,4)+f(3,1,5,2,8,10,3,4)+
f(5,1,5,2,9,6,7,4)+f(2,1,5,2,6,8,7,4)+f(9,1,5,9,6,8,7,4)+
f(8,1,2,9,6,8,7,4)+f(1,1,2,9,8,10,7,4)+f(5,5,2,9,8,10,7,4)+
f(6,2,9,6,8,10,7,4)+f(10,1,2,9,6,10,7,4)+f(10,1,5,6,8,10,7,4)+
f(1,1,5,2,9,3,7,4)+f(6,5,2,9,6,3,7,4)+f(8,5,2,9,8,3,7,4)+
f(5,1,5,6,10,3,7,4)+f(2,1,2,6,10,3,7,4)+f(9,1,9,6,10,3,7,4)+
f(8,1,6,8,10,3,7,4)+f(10,5,2,9,10,3,7,4)+f(3,1,2,9,6,3,7,4)+
f(3,1,5,6,8,3,7,4)+f(3,2,9,8,10,3,7,4)+f(7,1,5,2,8,10,7,4)+
f(7,1,5,9,8,3,7,4)+f(7,5,9,6,10,3,7,4)+f(7,2,6,8,10,3,7,4)+
f(4,1,5,2,9,6,10,4)+f(4,1,2,9,6,8,3,4)+f(4,5,9,6,8,10,7,4)+
f(4,5,2,6,8,3,7,4)+f(4,1,5,2,10,3,7,4)+f(4,1,9,8,10,3,7,4).$

The Igusa invariants are
$$I_2 = \frac{h_{12}}{h_{10}},$$
$$I_4 = h_4,$$
$$I_6 = \frac{h_{16}}{h_{10}},$$
$$I_{10} = h_{10},$$
and the absolute Igusa invariants are 
$$j_1 =\frac{I_2^5}{I_{10}},$$
$$j_2= \frac{I_4 I_2^3}{I_{10}},$$
$$j_3= \frac{I_6 I_2^2}{I_{10}}.$$


\begin{thebibliography}{999999}



\bibitem[CL01]{CL} Cohn, Henry; Lauter, Kristin.  {\it Generating Genus 2 Curves with 
Complex Multiplication}, Internal Technical Report, January 2001. 

\bibitem[GZ84]{GZ} Gross, Benedict H.; Zagier, Don B. {\it On singular moduli}. 
J. Reine Angew. Math. {\bf 355} (1985), 191--220. 

\bibitem[Sh94]{Shimura} Shimura, Goro. Introduction to the Arithmetic Theory of Automorphic 
Functions, Princeton University Press, 1994.

\bibitem[Sp94]{Spallek} Spallek, Anne-Monika.  {\it Kurven vom Geschlecht 2 und ihre Anwendung in Public-Key-Kryptosystemen}, Preprint No. 18 (1994), Institut f\"ur Experimentelle Mathematik, Essen.

\bibitem[vW99]{vW}  van Wamelen, Paul. {\it Examples of genus two CM curves 
defined over the rationals}. Math. Comp. {\bf 68} (1999), no. 225, 307--320.

\bibitem[We03]{Weng}  Weng, Annegret. {\it Constructing hyperelliptic curves of genus 2 
suitable for cryptography}. Math. Comp. {\bf 72} (2003), 435-458.


\end{thebibliography}
\end{document}